\documentclass[11pt]{article}
\usepackage{mathrsfs,amsfonts,amsmath}
\usepackage{bbm}
\RequirePackage{ifthen,calc}
\usepackage{theorem}
\usepackage[dvips]{color}
\usepackage{graphicx}
\usepackage{graphicx}
\usepackage{mathptmx}      
\usepackage{xcolor}%
\usepackage{ulem}%
\usepackage{graphicx,amsmath,amsfonts,amssymb,color,bbm}%

\def\R{{\mathbb R}}

\def\E{{\mathbb E}}
\def\P{{\mathbb P}}

\catcode`@=12

\newtheorem{thm}{\noindent Theorem}[section]
\newtheorem{lem}{\noindent Lemma}[section]
\newtheorem{cor}{\noindent Corollary}[section]

\newtheorem{con}{\noindent  Condition}[section]

\newtheorem{remark}{Remark}[section]
\newtheorem{definition}{Definition}[section]

\theoremheaderfont{\normalfont\bfseries}
\theorembodyfont{\slshape} {\theorembodyfont{\rmfamily}} {\theorembodyfont{\rmfamily}
\newtheorem{defn}{\noindent Definition}[section]}
{\theorembodyfont{\rmfamily} } {\theorembodyfont{\rmfamily}
}
{\theorembodyfont{\rmfamily} } {\theorembodyfont{\rmfamily}
\newtheorem{rem}{\noindent Remark}[section]}
{\theorembodyfont{\rmfamily} } {\theorembodyfont{\rmfamily} } {\theorembodyfont{\rmfamily}
}

 \def\beqlb{\begin{eqnarray}}\def\eeqlb{\end{eqnarray}}
 \def\beqnn{\begin{eqnarray*}}\def\eeqnn{\end{eqnarray*}}
 
 \numberwithin{equation}{section}
\def\qed{\hfill$\square$\smallskip}

 \newcommand{\pend}{\hfill \thicklines \framebox(5.5,5.5)[l]{}}
 \newenvironment{proof}{\noindent {\sc  Proof} \rm}{\pend}

\begin{document}
\title{Multidimensional Sticky Brownian Motions:  Tail Behaviour of the Joint Stationary Distribution }
\author{Hongshuai Dai$^{a,b}$ and Yiqiang Q. Zhao$^b$
 \\ {\small a. School of Statistics, Shandong University of Finance and Economics, Jinan, Shandong
250014, P.R. China}\\
{\small b. School of Mathematics and Statistics,
Carleton University, Ottawa, ON, Canada K1S 5B6}}
\maketitle
\begin{abstract}
Sticky Brownian motions, as time-changed semimartingale reflecting Brownian motions, have various applications in many fields, including queuing theory and mathematical finance. In this paper, we are concerned about the stationary distributions of a multidimensional sticky Brownian motion, provided it is stable.  We will study the large deviations principle for stationary distribution and the tail behaviour of the joint stationary distribution.
\end{abstract}

\small {{\bf MSC(2000):}  60K25, 60J10\\[1mm]


\section{Introduction}



Since the works of Feller \cite{F1952,F1954,F1957}, sticky  Brownian motions  have been  explored extensively, for example, see It\^{o} and McKean \cite{IM1963, IM1965}. Recall  that  a sticky Brownian motion on the half-line is the process evolving as a standard Brownian motion away from zero and reflecting at zero after spending a random amount of time there.  It is known that a sticky Brownian motion arises as a time change of a semimartingale reflecting Brownian motion (SRBM), which reflects at zero instantaneously, and  it describes the scaling limit of random walks on the natural numbers whose jump rate at zero is significantly smaller than that at positive sites. From a queuing theory perspective, this kind of process is quite an interesting one with many applications.  Welch~\cite{W1964} introduced an exceptional service for the first customer in each busy period and showed that a sticky Brownian motion on the half-line can be a heavy traffic limit. Later, with different exceptional service mechanisms, the same heavy traffic limit, or the sticky Brownian motion, was confirmed for other single server queuing models by Lemoine \cite{L1975},  Harrison and Lemoine \cite{HL1981}, Yamada \cite{Y1994}, and Yeo \cite{Y1961}.

 Recently, R\'acz and Shrocnikov \cite{RS2015}  introduced multidimensional sticky Brownian motions, which are a natural multidimensional  extension of sticky Brownian motions on the half-line.  As shown in \cite{RS2015}, a multidimensional sticky Brownian motion can also be written as a time-changed multidimensional semimartingale reflecting Brownian motion. Multidimensional sticky Brownian motions have many potential applications in both queuing theory and mathematical finance. For example, as pointed out by  R\'acz and Shrocnikov \cite{RS2015}, it can be used as an approximation of certain particle movement systems.

Stationary distributions of the multidimensional SRBM have attracted a lot of interest.  When $R$ is an $\mathcal{M}$-matrix and $R^{-1}\mu<0$,  Majewski \cite{M1998}, and Avram, Dai and Hasenbein \cite{ADH2001} established the large deviations principle (LDP).  The $\mathcal{M}$-matrix condition can be relaxed, for example, Dupuis and Ramanan \cite{DR2002}, under more general conditions,  studied a time-reversed representation for the tail probabilities of an SRBM.  We are interested in the tail behaviour of stationary distributions. Many efforts have been made to study this topic and some results have been obtained, most of which are related to special multidimensional cases, including the two-dimensional case and the skew symmetric case. Intuitively, we can discuss them by the large deviations principle. Then, the problem is reduced to finding the rate function, which is formulated as a variational problem. However, it is known that, in general, with the exception of some special cases, it is very difficult to analytically solve these variational problems. Some discussions about why higher dimensional $(\geq 3)$ cases are difficult have been carried out, for example, see Avram, Dai and Hasenbein  \cite{ADH2001}.  Hence, additional work is needed to study this problem for the multidimensional SRBM. For the two-dimensional SRBM,  Dai and Miyazawa \cite{DM2011}, Franceschi and Raschel \cite{FR2017},  Dai, Dawson and Zhao~\cite{DDZ}, and Franceschi and Kurkova \cite{FK2016} studied the tail asymptotics of the marginal distributions of the SRBM, and obtained the decay rate of the marginal distributions. At the same time, being inspired by the two-dimensional case and some special multidimensional cases, we note that some conjectures on the tail properties of the stationary marginal distributions of the multidimensional SRBM have been discussed. Miyazwa and Kobayashi \cite{MK2011} conjectured on the decay rate of the marginal distribution in an arbitrary direction of the multidimensional SRBM.  Motivated by the above arguments, in this work, we also study some of the tail properties of the stationary distributions of the multidimensional sticky Brownian motion.

In a recent paper, Dai and Zhao \cite{DZ2018} obtained exact tail asymptotics and asymptotic independence of a two-dimensional sticky Brownian motion. The main tools applied in \cite{DZ2018} were the kernel method, extreme value theory and copula.  In this paper, we extend those ideas for the general multidimensional case.  We first discuss the LDP for the sticky Brownian motion, and then study the tail behaviour of the joint stationary distribution of the process.

The rest of this paper is organized as follows: in Section~\ref{sec:2}, we first state some preliminaries related to multidimensional sticky Brownian motions, and then we study some basic properties of the stationary distributions of multidimensional sticky Brownian motions. In Section~\ref{sec:3}, we establish the LDP for the sticky Brownian motion. In Section~\ref{sec:4}, we study the stationary behaviour for the joint stationary distribution.  In Section~\ref{sec:5},  we discuss  exact tail asymptotics  for some special cases of the multidimensional sticky Brownian motion.

\section{Sticky Brownian motion} \label{sec:2}

In this section, we introduce some preliminaries related to multidimensional sticky Brownian motions. We first recall the definition  of the semimartingale reflecting Brownian motion (SRBM). SRBM models arise as an approximation for queuing networks of various kinds (see, for example, Williams~\cite{W1995,W}). A $d$-dimensional SRBM, denoted as $\tilde{Z}=\{\tilde{Z}(t),t\geq 0 \}$, is defined as follows:
 \beqlb\label{def1}
 \tilde{Z}(t)=X(t)+R L(t), \;\textrm{for}\;t\geq 0,
 \eeqlb
where  $\tilde{Z}(0)=X(0)\in\R^d_+$,  $X$ is an unconstrained Brownian motion with drift vector
$\mu=(\mu_1,\cdots,\mu_d)'$ and covariance matrix $\Sigma=(\Sigma_{i,j})_{d\times d}$, $R=(r_{ij})_{d\times d}$ is a $d \times d$ matrix specifying the reflection behaviour at the boundaries,  and
$L=\{L(t)\}$ is a $d$-dimensional process with the local times $L_1, \ldots,
L_d$ such that:
\begin{itemize}
\item[(i)]the local time $L_i$ is continuous and non-decreasing with $L_i(0)=0$;
\item[(ii)] $L_j$ only increases at times $t$ for which
$\tilde{Z}_i(t)=0$, $i=1,\ldots, d$;
\item[(iii)] $\tilde{Z}(t)\in\R^d_+$, $t\geq 0$.
\end{itemize}

The existence of an SRBM has been studied extensively, for example, Taylor and Williams~\cite{TW1993}, and
Reiman and Williams~\cite{RW1988}.  It was proved in \cite{RW1988,TW1993} that, for a given set of
data $(\Sigma, \mu,R)$, with $\Sigma$ being positive definite,
there exists an SRBM for each initial distribution of $\tilde{Z}(0)$, if
and only if, $R$ is completely $\mathbb{S}$ (see, for example, Taylor and Williams~\cite{TW1993} for the definitions of matrix classes).  Furthermore, when $R$
is completely $\mathbb{S}$, the SRBM is unique in distribution for each given initial distribution. It is well-known that a necessary condition (see, for example, Harrison and Williams~\cite{HW1987}, or Harrison and Hasenbein~\cite{HH2009}) for the
existence of the stationary distribution for $\tilde{Z}$ is
\beqlb\label{M-1} R \;\textrm{is non-singular and }\;R^{-1}\mu<0.
\eeqlb

We note that an SRBM does not spend time on the boundary.  Conversely, a sticky Brownian motion would spend a duration of  time on the boundary.  For the one-dimensional case, Feller~\cite{F1952,F1954,F1957} first observed the sticky boundary behaviour for diffusion processes and studied the problem that describes domains of the infinitesimal generators associated with a strong Markov process $\tilde{X}$ in $[0,\,\infty)$.  Moreover, $\tilde{X}$ behaves like a standard Brownian motion  in $(0,\;\infty)$, while at $0$, a possible boundary behaviour is described by
\beqlb\label{2-a1}
f'(0+)=\frac{1}{2u}f''(0+),
\eeqlb
where $u\in(0,\infty)$ is a given and fixed constant and $f$ are functions belonging to the domain of the infinitesimal generator of $\tilde{X}$.
The second derivative $f''(0+)$ measures the ``stickiness" of $\tilde{X}$ at $0$. For this reason, the process $\tilde{X}$ is called a sticky Brownian motion, which is also referred to as a sticky reflecting Brownian motion in the literature. It\^{o} and Mckean~\cite{IM1963} first constructed the sample paths of $\tilde{X}$.  They showed that $\tilde{X}$ can be obtained from a one-dimensional SRBM $\tilde{Z}$ by the time-change $t\to T(t):=S^{-1}(t)$,  where $S(s)=s+\frac{1}{u}L_s$ for $s>0$, or $T(t)=s$ is determined by the equation $t = s+\frac{1}{u}L_s$.
For more information about sticky Brownian motions on the half-line, refer to Engelbert and Peskir~\cite{EP2014} and the references therein.

 R\'acz and Shrocnikov~\cite{RS2015} introduced multidimensional sticky Brownian motions and proved the existence and uniqueness of the multidimensional sticky Brownian motion.  Similar to a sticky Brownian motion on the half-line,  let
\beqlb\label{2-a3}
S(t)=t+\sum_{i=1}^du_iL_i(t),
\eeqlb
where $u_i\in(0,\;\infty)$, $i=1,\ldots,d$, are given and fixed constants. For convenience, let $u=(u_1,\cdots,u_d)'$.  Let $T(\cdot)$ be the inverse of $S(t)$, that is,
\beqlb\label{T} T(t)=S^{-1}(t).\eeqlb Then, it  follows from Kobayashi \cite[Lemma 2.7]{K2011} and the equation \eqref{2-a3} that  $T$ has continuous paths and $\lim_{t\to\infty}T(t)=\infty.$  Furthermore,  $0<T(1)\leq 1$. Then, a multidimensional sticky Brownian motion can be defined as:
\beqlb\label{1-4}
Z(t)=\tilde{Z}\big(T(t)\big).
\eeqlb

This type of process finds applications in the fields of queuing theory and mathematical finance. In the queuing field,  it is well known that the SRBM is a heavy traffic limit for many queuing networks such as open queuing networks. As discussed in the introduction, in the setting for single server queues, a sticky Brownian motion on the half-line can be served as a heavy traffic limit of a queuing system with exceptional service mechanisms. It is reasonable to expect that a multidimensional sticky Brownian motion serves as a heavy traffic limit for such multidimensional queuing networks with appropriately defined exceptional service mechanisms.

In the rest of this section, we study some of the properties of the stationary distributions of multidimensional sticky Brownian motions.  We first establish the so-called basic adjoint relation (BAR), which establishes some connections between  the joint stationary distribution and the boundary stationary measures defined below.  In particular, in the two-dimensional case, the BAR can be used to study exact tail asymptotics for the marginal stationary distributions and the boundary stationary measures of a sticky Brownian motion (see, Dai and Zhao \cite{DZ2018}).

In the rest of this paper, we assume that $Z(0)$ follows the stationary distribution $\pi$ of $\{Z(t)\}$. Furthermore, for the stationary measure $\pi$, we define the moment generating function (MGF) $\Phi(\theta)$ by
\beqnn
    \Phi(\theta)=\int_{\R_+^d}\exp\{<\theta,x>\} \pi(dx).
\eeqnn
Similar to the SRBM, $\Phi(\theta)$ is closely related to the MGFs of various boundary measures, which are defined below.
For any set $A\in\mathscr{B}(\R_+^d)$, define
\beqlb\label{2-25}
V_i(A)&&=\E_{\pi}\bigg[\int_0^{T(1)} 1_{\{\tilde{Z}(s)\in A\}}dL_i(s)\bigg].
\eeqlb
At the same time, for any Borel  set $B\in\mathscr{B}(\R_+^d)$, we define the joint measure for the time-change as:
\beqnn
    V_0(B)=\E_{\pi}\bigg[\int_0^1 1_{\{Z(s)\in B\}}dT(s)\bigg].
\eeqnn
According to Lemma~\ref{thm-1a} below, all $V_i$, $i=0, \cdots,d$, are finite measures on $\R_+^d$.
Then, we can define MGFs $\Phi_i(\theta)$ for $V_i$, $i=0,1,\cdots,d$, by
\beqnn
\Phi_i(\theta)=\int_{\R_+^d}\exp\{<\theta,x>\}V_i(dx).
\eeqnn
For these measures,  we have the following BAR:
\begin{lem}\label{thm-1a}\quad
\begin{itemize}
\item[(1)]The boundary measures $V_i$, $i=1,\cdots,d$, and the joint measure $V_0$ for the time-change are all finite.
\item[(2)]The MGFs of $V_i$, $i=0,1,\cdots,d$, have the following BAR: for any ${\bf\theta}\in \R_{-}^d = \{ \theta = (\theta_1,\cdots, \theta_d)': \theta_i\leq 0\}$,
\beqlb\label{2-5-a}
-\Psi_{X}({\bf\theta})\Phi_0({\bf\theta})=\sum_{i=1}^d\Phi_i({\bf\theta})<{\bf\theta}, R_i>,
\eeqlb
where $R_i$ is the $i$th column of the reflection matrix $R$, and $\Psi_{X}(\theta)$ is the L\'evy exponent of the multidimensional Brownian vector $X(1)$, i.e.,
\beqnn
\Psi_{X}(\theta)=<\theta,\mu>+\frac{1}{2}<\theta,\Sigma\theta>.
\eeqnn
\end{itemize}
\end{lem}
\noindent\underline{\proof:} Since $Z(0)$ follows the stationary distribution $\pi$, for any $t\in\R_+$,
\beqnn\label{2-27}
\P(Z(t)\leq z)=\P(Z\leq z).
\eeqnn
We note that $\{Z(t)\}$ is a semimartingale. Since $T(t)$ is continuous and $S(t)$ is strictly increasing,  it follows from Kobayashi~\cite[Corollary 3.4]{K2011} that  if $f:\R^d\to \R$ is $C_{b}^2$ function, then
\beqlb\label{3-5}
&&f(Z(t))-f(Z(0))=\sum_{i=1}^d \mu_i\int_0^{T(t)}\frac{\partial f}{\partial x_i}(\tilde{Z}(u))du+ \sum_{i,j=1}^d \int_0^{T(t)}r_{ji}\frac{\partial f}{\partial x_j}\big(\tilde{Z}(u)\big)dL_i(u)\nonumber
\\&&\hspace{1cm}+\sum_{i=1}^d\int_0^{T(t)}\frac{\partial f}{\partial x_i}\big(\tilde{Z}(u)\big)dX_i(u)+\frac{1}{2}\sum_{i,j=1}^d\Sigma_{i,j}\int_0^{T(t)}\frac{\partial^2 f}{\partial x_i\partial x_j}\big(\tilde{Z}(u)\big)du.
\eeqlb
Hence, we have
\beqlb\label{2-6}
&&\sum_{i=1}^d \mu_i\E_{\pi}\bigg[\int_0^{T(t)}\frac{\partial f}{\partial x_i}(\tilde{Z}(u))du\bigg]+ \sum_{i,j=1}^d \E_\pi\bigg[\int_0^{T(t)}r_{ji}\frac{\partial f}{\partial x_j}\big(\tilde{Z}(u)\big)dL_i(u)\bigg]\nonumber
\\&&\hspace{1cm}+\frac{1}{2}\sum_{i,j=1}^d\Sigma_{i,j}\E_\pi\bigg[\int_0^{T(t)}\frac{\partial^2 f}{\partial x_i\partial x_j}\big(\tilde{Z}(u)\big)du\bigg]=0.
\eeqlb
Next, we prove the first part of this lemma.   From \eqref{2-25}, we get that for all $i=1,\cdots,d$,
\beqlb
V_i(\R_+^d)=\E_\pi\Big[ L_i\big(T(1)\big)\Big],
\eeqlb
and
\beqnn
V_0(\R_+^d)=\E_\pi\Big[T(1)\Big].
\eeqnn
Hence, it suffices to prove that for any $i\in\{1,\cdots,d\}$,
\beqlb\label{2-28}
\E_\pi\Big[ L_i\big(T(1)\big)\Big]<\infty,
\eeqlb
and
\beqlb\label{2-28a}
\E_\pi\Big[T(1)\Big]<\infty.
\eeqlb
Since $T(1)\leq 1$, in order to prove \eqref{2-28}, we only need to show that
\beqlb\label{2-28b}
\E_\pi\Big[ L_i\big(1\big)\Big]<\infty.
\eeqlb
It follows from Dai and Harrison \cite[Proposition 3]{DH1992} that \eqref{2-28b} holds. Then \eqref{2-28} follows.  At the same time, from the relationship between $T(\cdot)$ and $S(\cdot)$, we get that
\beqlb\label{2-a2}
T(t)=t-\sum_{i=1}^du_iL_i(T(t)).
\eeqlb
Hence,
\beqlb\label{2-aa2}
\E_{\pi}[T(1)]=1-\sum_{i=1}^d\E_{\pi}u_iL_i(T(1)).
\eeqlb
Combining \eqref{2-28} and \eqref{2-aa2} leads to \eqref{2-28a}.
\newline

Taking $f(x_1,\cdots,x_d)=\exp\{\sum_{i=1}^d\theta_i x_i\}$ with $\theta_i\leq 0$, $i=1,\cdots,d$, in the equation \eqref{2-6} can prove the second part of this lemma. \qed

\begin{remark}\label{2-rem} Let $ C^{2}_{b}(\R^{d}_{+})$  be the set of functions $f$ on $\R_+^d$  such that $f$, its first order derivatives, and its second order derivatives are
bounded and continuous. For any $f\in C^{2}_{b}(\R^{d}_{+})$, it follows from \eqref{2-6} that
\beqnn
\int_{\R_+^d}\mathcal{L}f(x) V_0(dx)+\sum_{i=1}^2\int_{\R_+^d}<\bigtriangledown f(x), R_{i}>V_i(dx)=0,
\eeqnn
where
\beqnn
\mathcal{L}f(x)=\frac{1}{2}\sum_{i,j=1}^d\Sigma_{i,j}\frac{\partial^2f}{\partial x_i\partial x_j}(x)+\sum_{j=1}^d\mu_j\frac{\partial f}{\partial x_j}(x),
\eeqnn
and $\bigtriangledown f(x)$ is the gradient of $f$.
From Dai and Kurtz~\cite[Theorem 1.4]{DK2003} (or Braverman, Dai and Miyazawa \cite[Lemma 2.1] {BDM2017}), we can get that $V_0(\cdot)/\E_\pi\big[T(1)\big]$, and $V_i(\cdot)/\E_\pi\big[L_i(T(1))\big]$, $i=1,\cdots,d$, are the stationary distribution, and the boundary distributions of the corresponding reflecting Brownian motion $\tilde{Z}$, respectively.
\end{remark}

The following corollary immediately follows from the proof to Lemma~\ref{thm-1}.
\begin{cor}\label{cor-2}   $\E\big[L_i(T(1))\big]$, $i=1,\cdots,d$, satisfy
\beqnn
\Big[\mu u'-R^T\Big]\tilde{L}=\mu,
\eeqnn
where
\beqnn
\tilde{L}=\bigg(\E\Big(L_1\big(T(1)\big)\Big),\cdots,\E\Big(L_d\big(T(1)\big)\Big)\bigg)'.
\eeqnn

\end{cor}
\underline{\proof:}Let $f(x_1,\cdots,x_d)=\exp\{\theta x_1\}$ with $\theta<0$ and $x_1\geq 0$. Then we have that
\begin{equation}\label{2-32}
f'_i(x_1,\cdots,x_d)=\left\{
             \begin{array}{lr}
             \theta\exp\{\theta x_1\}, \;\textrm{if}\; i=1&
               \\  0, \hspace{1.5cm}\textrm{other}&
             \end{array}
\right.
\end{equation}
and
\begin{equation}\label{2-32a}
f''_{i,j}(x_1,\cdots,x_d)=\left\{
             \begin{array}{lr}
             \theta^2\exp\{\theta x_1\}, \;\textrm{if}\; i=j&
               \\  0, \hspace{1.5cm}\textrm{other}.
                \end{array}
\right.
\end{equation}
Hence, combining  \eqref{2-6},  \eqref{2-32} and \eqref{2-32a} gives
\beqlb\label{2-31}
-\Psi_{X}(\theta,0,\cdots,0)\Phi_0(\theta,0,\cdots,0)= \sum_{i=1}^d\Phi_i(\theta,0,\cdots,0)r_{i1}\theta.
\eeqlb
Dividing $\theta<0$ at both sides of \eqref{2-31} and letting $\theta\to 0$, we get
\beqlb\label{2-33}
-\mu_1\E_{\pi}[T(1)]=\sum_{i=1}^d\E_\pi[L_i(T(1))]r_{i1}.
\eeqlb
Symmetrically, let $f(x_1,\cdots,x_d)=\exp\{\theta_i x_i\}$ with $\theta_i<0$ and $x_i\geq 0$, $i=2,\cdots,d$. Similar to \eqref{2-33}, we get
\beqlb\label{2-34}
-\mu_i\E_{\pi}[T(1)]=\sum_{j=1}^d\E_\pi[L_j(T(1))]r_{ji}.
\eeqlb
Then, combining \eqref{2-33} and \eqref{2-34} yields
\beqlb\label{2-34a}
\Big[\mu u'-R^T\Big]\tilde{L}=\mu.
\eeqlb
The lemma is proved.
\qed
\newline

Below, we state the main result of this section. The sticky Brownian motion defined by \eqref{1-4} satisfies the following BAR:
\begin{thm}\label{thm-1}
\beqlb\label{2-5}
-\Psi_{X}({\bf\theta})\Phi({\bf\theta})=\sum_{i=1}^d\Phi_i({\bf\theta})\bigg(<{\bf\theta}, R_i>-u_i\Psi_{X}({\bf\theta})\bigg).
\eeqlb
\end{thm}

\noindent\underline{\proof:}
For any Borel set $B\in\mathscr{B}(\R_+^d)$, we have
\beqlb\label{2-8}
\pi(B)&=&\E_\pi\Big[\int_0^1 \mathbbm{1}_{\{Z(s)\in B\}}ds\Big]=\E_\pi\Big[\int_0^{T(1)} \mathbbm{1}_{\{\tilde{Z}(s)\in B \}}dS(s)\Big]\nonumber
\\&=& \E_\pi\Big[\int_0^{T(1)} \mathbbm{1}_{\{\tilde{Z}(s)\in B \}}dt\Big]+\sum_{i=1}^du_i\E_\pi\Big[\int_0^{T(1)} \mathbbm{1}_{\{\tilde{Z}(s)\in B\}}dL_i(s)\Big]\nonumber
\\&=& V_0(B)+\sum_{i=1}^d u_iV_i(B).
\eeqlb

From \eqref{2-5-a} and \eqref{2-8}, we can get \eqref{2-5}. The proof is completed. \qed
\begin{remark}
The BAR plays an important role in analyzing the asymptotic properties of various stationary distributions.
For the two-dimensional case, based on the BAR, Dai and Zhao \cite{DZ2018} applied the kernel method to get exact tail behaviours of the marginal stationary distributions and  various boundary measures.  On the other hand, we note that when all $u_i=0$, for $i=1,\cdots,d,$ the  sticky Brownian motion $Z$ reduces to  an SRBM.  In this case, the BAR still holds, which has been obtained (see,  for example, Harrison and Williams \cite{HW1987}).
\end{remark}

\section{Large Deviations for Sticky Brownian Motion}
\label{sec:3}

In this section, we study the large deviations principle for the stationary  distribution $\pi$ of $Z$. We first recall the definition of $LDP$, see, for example, Varadhan \cite[Defintion 2.1]{V1984}.
\begin{defn}
A sequence of probability measures $\{\mu_n\}$ defined on a complete separable metric space $(\mathscr{H},\mathscr{B})$ is said to satisfy the LDP with speed $\{\xi_k\}$ and rate function $J$ if, for all $\Theta\in\mathscr{B}$ and $\lim_{k\to\infty}\xi_k=\infty$,
\beqnn
\limsup_{n\to\infty}\frac{1}{\xi_n}\log \mu_n(\Theta)\leq -\inf_{x\in\bar{\Theta}} J(x),
\eeqnn
and
\beqnn
\liminf_{n\to\infty}\frac{1}{\xi_n}\log \mu_n(\Theta)\geq -\inf_{x\in \Theta^0} J(x),
\eeqnn
where $J:\mathscr{H}\to[0,\infty]$ is a function with compact level sets, and $\bar{\Theta}$ (respectively $\Theta^0$) is the closure (respectively interior) of $\Theta$.  A sequence of random variables $\{X_n\}$ defined on some measure space taking values in a complete separable  metric space $(\mathscr{H},\mathscr{B})$  is said to satisfy an LDP, with rate function $J(\cdot)$, if the corresponding induced measures satisfy a LDP with the same rate function.
\end{defn}
To reach our objective, we need the contraction principle (see, for example, Amir and Ofer \cite[Theorem 4.2.1]{AO1993}). Here we briefly recall it.  Suppose that a sequence of random variables $\{X_n\}$ satisfies a large deviations principle with speed $\{\xi_k\}$  and good rate function $J$ in the topology $\mathscr{B}$, and $f:\mathscr{H}\to\mathscr{H}'$  is a continuous and measurable mapping to the topological space $(\mathscr{H}',\mathscr{B}')$. Then the contraction principle states that the sequence $\{f(X_n)\}$ satisfies a large deviation principle with speed $\{\xi_k\}$  and good rate function $J':\mathscr{H}'\to[0,\infty]$ given for $x'\in\mathscr{H}'$ by
\beqnn
J'(x')=:\inf_{x\in E,x'=f(x)} J(x),
\eeqnn
in the topology $\mathscr{B}'$.

Here, we need the LDP for the SRBM. Let $\tilde{\mu}_n(B)=\tilde{\pi}(nB)$, where $\tilde{\pi}$ is the stationary distribution of the SRBM, and $\mathscr{A}([0,\infty);\R^d)$ be the corresponding sets of absolutely continuous functions on $[0,\infty)$ taking values in $\R^d$.  We note that the LDP for an SRBM has been studied extensively (see, for example,  Avram, Dai and Hasenbein \cite{ADH2001}, Majewski \cite{M1998}, Dupuis and Ramanan \cite{DR2002} and the references therein ).  However, there is no LDP established in other literature when $\Sigma$ is  completely-$\mathbb{S}$.  In this paper, we  also study the LDP for $Z$ under some mild conditions.  We will study the LDP under the conditions  in  Dupuis and Ramanan \cite{DR2002}. We first recall these conditions.
\begin{con}\label{SRBM-C1}
$R$ is invertible and the associated Skorokhod Map $\Gamma$ is Lipschitz continuous (with respect to the topology of uniform convergence on compact sets) and is defined for every  $\Psi\in \mathcal{C}_+([0,\infty):\R^d)$.
\end{con}

\begin{rem}
To satisfy Condition \ref{SRBM-C1}, $R$ must, in particular, be completely-$\mathbb{S}$.  Discussions about general assumptions that ensure Condition \ref{SRBM-C1}  can be found in Dupuis and Ramanan \cite{DR1999a, DR1999b, DR2002} and Dupuis and Ishill \cite{DI1991}.
\end{rem}

\begin{con}\label{SRBM-C2}
Define $\mathscr{L}=\big\{-\sum \alpha_iR_i: \alpha_i\geq 0\big\}$, where $R_i$ is the $i$th column of the matrix $R$.  Assume that $\mu\in \mathscr{L}$.
\end{con}

\begin{rem}
If $R$ is invertible, then Condition \ref{SRBM-C2} is equivalent to the inequality \eqref{M-1}.
\end{rem}

 The following lemma comes from Dupuis and Ramanan \cite{DR2002}.
 \begin{lem} \label{LDP-SRBM} Assume that the SRBM is  such that $\Sigma$ is positive definite and Conditions \ref{SRBM-C1} and \ref{SRBM-C2} are satisfied. Le $\Gamma$ be the associated Skorohod map.  Then, $\{\tilde{\mu}_n\}$ satisfies the LDP with speed function $n$ and the rate function $\tilde{V}(x)$ that is given by
 \beqnn
\tilde{V}(x)=\inf_{\phi\in\mathscr{A}([0,\infty):\R^d): \phi(0)=0,\phi\in\Gamma(\Psi):\tau_x <\infty} \int_0^{\tau_x} L(\dot{\Psi}(s))ds,
 \eeqnn
where
\beqnn
L(\beta)=\frac{1}{2}(\beta-b)'\Sigma^{-1}(\beta-b),
\eeqnn
and
\beqnn
\tau_x=\inf\{t\geq0:\phi(t)=x\},
\eeqnn
and where $\Gamma(\psi)$  is the set of images of $\psi$ under the Skorohod Map (SM) that is associated with $\Sigma$.
 \end{lem}
Next, we state the LDP for the stationary distributions $\pi$ of the sticky Brownian motion $Z$. Let $\mu_n(B)=\pi(nB)$.  It follows from equation \eqref{2-a3} that $S(t)$ is strictly increasing.  Then it follows from Lemma 2.7 in \cite{K2011} that $T$ has continuous sample paths.    Hence, let  $\mathbb{T}:C\big([0,\infty),\R^d\big)\to C\big([0,\infty),\R^d\big)$ be a continuous  function such that, for any $\omega$,
\beqnn \mathbb{T}: &C\big([0,\infty),\R^d\big)&\to C\big([0,\infty),\R^d\big)  \\
&\tilde{Z}(t)&\to Z(t)=\tilde{Z}(T(t)).
\eeqnn
Hence, from the contraction principle and Lemma 3.1, we have the following LDP for $\mu_n$.
 \begin{thm} \label{LDP} Assume that $\Sigma$ is positive definite and Conditions \ref{SRBM-C1} and \ref{SRBM-C2}  are satisfied. Then $\{\mu_n\}$ satisfies the LDP with the rate function $V(x)$, given by
 \beqnn
 V(x)=\inf_{x'\in C\big([0,\infty),\R^d\big),x=\mathbb{T}(x')} \tilde{V}(x').
 \eeqnn
 \end{thm}

\begin{rem}\label{LDP-R}
From Theorem \ref{LDP}, we can see that for any measurable set $B\subset \R^d$,
\beqnn
&&\limsup_{n\to\infty}\frac{1}{n}\log \P\big(Z\in nB\big)\leq \alpha_{\bar{B}};\\
&&\liminf_{n\to\infty}\frac{1}{n}\log \P\big(Z\in nB\big)\geq \alpha_{B^o},
\eeqnn
where $ \alpha_{\bar{B}}=-\inf_{x\in \bar{B}}V(x)$ and $\alpha_{B^o}=-\inf_{x\in B^0}V(x)$.
 \end{rem}
 To study the tail behaviour of the joint stationary distribution, we need the tail properties of the marginal $Z_i$, $i\in\{1,\cdots,d\}$.  From Theorem \ref{LDP} and Remark \ref{LDP-R}, we have the following corollary.
 \begin{cor}\label{con-1}
Assume that $\Sigma$ is positive definite and Conditions \ref{SRBM-C1} and \ref{SRBM-C2}  are satisfied. Then for any $i\in\{1,\cdots,d\}$,
\beqlb\label{4-2}
-\lim_{x\to\infty}\frac{1}{x}\log\P\{Z_i\geq x\}=\alpha_i.
\eeqlb
 \end{cor}

\section{Tail Behaviour of the Joint Distribution} \label{sec:5}

It is well known that if we have a multivariate Gaussian vector, where the correlation coefficients are strictly less than $1$, then it is asymptotically independent (see Definition \ref{defn-AI} below). On the other hand, we note that in the interior of the first quadrant $\R_+^d$, the sticky Brownian motion $Z$ behaves like the Brownian motion. Hence, it is expected that, under some mild conditions, $Z$ is also asymptotically independent.  In this section, we discuss the asymptotic independence of $Z$. In the rest of this paper, we first  assume that all the correlation coefficients $\rho_{X_iX_j}<1$, $i,j\in\{1,\cdots,d\}$ where  $X(1)=(X_1,\cdots,X_d)'$.

To study the tail behaviour of the joint stationary distribution, we mainly use the copula. 
For any multidimensional distribution $\tilde{F}$ with marginal distributions $\tilde{F}_i$, $i=1,\cdots,d$, the copula associated with $\tilde{F}$ is a distribution function
$C:[0,\;1]^d\to[0,\;1]$ satisfying
\beqnn
\tilde{F}({\bf x})=C\big(\tilde{F}_1(x_1),\cdots,\tilde{F}_d(x_d)\big).
\eeqnn
For more information on copula, we refer the reader to Joe \cite{J1997}.
Therefore, if $C(\cdot)$ is a copula, then it is a multivariate distribution with all univariate marginal distributions being $U(0,\;1)$, or the joint distribution of a multivariate uniform random vector.  
It is also well known that for continuous multivariate distributions, the univariate margins and the multivariate or dependence structure can be separated, and the multivariate structure is represented by a copula.

Next, we discuss the tail behaviour of the joint stationary distribution. We first recall some definitions.

\begin{definition}[Domain of Attraction] Assume that $\big\{X_n=(X_1^{(n)},\ldots,X_d^{(n)})'\big\}$ are independent and identical distributed (i.i.d.) multivariate random vectors with common distribution $\tilde{F}(\cdot)$ and the marginal distributions $\tilde{F}_i(\cdot)$, $i=1,\ldots,d$. If there exist normalizing constants $a_n^{(i)}>0$ and $b_i^{(n)}\in\R$, $1\leq i\leq d$, $n\geq 1$ such that, as $n\to\infty$,
\beqnn
\P\Big\{\frac{M_i^{(n)}-b_i^{(n)}}{a_i^{(n)}}\leq x_i,1\leq i\leq d\Big\}&&=\tilde{F}^n\Big(a_1^{(n)}x_{1}+b_1^{(n)},\ldots,a_d^{(n)}x_{d}+b_d^{(n)} \Big)
\\&&\to G(x_1,\ldots,x_d),
\eeqnn
where $M_i^{(n)}=\bigvee_{k=1}^{n}X_i^{(k)}$ is the componentwise maxima, then we call the distribution function $G(\cdot)$ a multivariate extreme value distribution function, and
$F$ is in the domain of attraction of $G(\cdot)$. We denote this by $\tilde{F}\in D(G)$.
\end{definition}

\begin{definition}\label{defn-AI}[Asymptotic Independence] Assume that the extreme value distribution function $G(\cdot)$ has the marginal distributions $G_i(\cdot)$, $i=1,\ldots,d$. If
\beqnn
\tilde{F}^n\Big(a_1^{(n)}x_1+b_1^{(n)},\ldots,a_d^{(n)}x_d+b_d^{(n)} \Big)
\to G(x_1,\ldots,x_d) =\prod_{i=1}^{d}G_i(x_i),
\eeqnn
then we say that $\tilde{F}(\cdot)$ is asymptotically independent.
\end{definition}

In the rest of this section, under some mild assumptions, we study the tail  behaviour of the joint stationary distribution $F(\cdot)$ of $Z$. As is standard for L\'evy-driven queueing networks, we assume below that the reflection matrix $R=I-P^T$, where $P$ is a substochastic matrix with its spectral radius strictly less than $1$, and $A^T$ is the transpose of an square matrix $A$.  From Condition 2.2 in Dupuis and Rananan \cite{DR2002},  we know that $R$ satisfies conditions 3.1 and 3.2.

To study the tail  behaviour of the joint stationary distribution $F$,  we first need to study the extreme value distribution of the univariate marginal stationary distribution $F_i(\cdot)$.  We have the following technical lemma.

\begin{lem}\label{con-2}  For any $i\in\{1,\cdots,d\}$,
\beqlb\label{4-1}
F_i(x)\in D(G_1),
\eeqlb
where
\beqlb\label{4-6}
 G_1(x)=\exp\{-e^{-x}\}.
\eeqlb
\end{lem}

\noindent \underline{\proof of Lemma \ref{con-2}:}~It follows from \eqref{4-2} and Corollary \ref{con-1} that
\beqlb\label{4-3}
\alpha_i=\lim_{x\to\infty}\frac{1}{x}\bigg(-\log \Big(1-F_i(x)\Big)\bigg).
\eeqlb
 From \eqref{4-3}, we obtain
\beqlb\label{4-7}
1-F_i(x)= \exp\{-\alpha_i x- o(\alpha_i x)\}.
\eeqlb
For convenience, let  $g_i(x)=\exp\{-o(\alpha_i x)\}$ and
\beqlb\label{4-26}\tilde{g}_i(x)=o(\alpha_i x).\eeqlb  Hence,
\beqlb\label{4-8}
g_i(x)=\exp\{-\tilde{g}_i(x)\}.
\eeqlb
Noting \eqref{4-26}, we can furthermore assume that $\tilde{g}_i(x)$ is twice continuously differential.   In such case, \eqref{4-3} and \eqref{4-8} suggest that
\beqlb\label{4-4}
1-F_i(x)\sim g_i(x)\exp\{-\alpha_i x\},\;\textrm{as}\; x\to\infty.
\eeqlb

Noting that
\beqlb\label{4-13}
\lim_{x\to\infty}\frac{\tilde{g}_i(x)}{\alpha_i x}=0,
\eeqlb
we get
\beqlb\label{4-12}
\lim_{x\to\infty}\tilde{g}_i(x)=\left\{
             \begin{array}{lr}
              K, \;\textrm{where}\; K \;\textrm{is a fixed and finite constant,}\;&
               \\  \infty.&
             \end{array}
\right.
\eeqlb
If $\lim_{x\to\infty}\tilde{g}_i(x)=\infty,$ then from \eqref{4-13} and the L'H$\hat{o}$spital rule, we get, as $x\to\infty$,
\beqlb\label{4-9}
\lim_{x\to\infty}\frac{\tilde{g}_i(x)}{\alpha_ix}=\lim_{x\to\infty}\frac{\tilde{g}_i'(x)}{\alpha_i}=0.
\eeqlb
Hence, from \eqref{4-12} and \eqref{4-9}, we get
\beqlb\label{4-20}
\lim_{x\to\infty}\tilde{g}'_i(x)=0.
\eeqlb
Furthermore, from \eqref{4-20}, we obtain
\beqlb\label{4-21}
\lim_{x\to\infty}\tilde{g}''_i(x)=0.
\eeqlb
Finally, from \eqref{4-8},  \eqref{4-4}, \eqref{4-9}, \eqref{4-20} and \eqref{4-21}, we get
\beqlb\label{4-5}
\lim_{x\to\infty}\frac{F''_i(x)\big(1-F_i(x)\big)}{\Big(F'_i(x)\Big)^2}=-1.
\eeqlb
Hence, from \eqref{4-5} and Proposition~1.1 in Resinck~\cite[pp. 40]{R1987}, we conclude that $F_i\in D(G_1)$. \qed

\begin{lem}\label{4-thm} For the sticky Brownian motion $Z=(Z_1,\cdots,Z_d)'$ with the stationary distribution function $F$,
\beqnn
F^n(a_i^{(n)} x_i+b_i^{(n)},i=1,\cdots,d )\to \Pi_{i=1}^dG_1(x_i),\;\text{as}\; n\to\infty,
\eeqnn
where $a^{(n)}_{i}$ and $b_i^{(n)}$ are normalizing constants.
\end{lem}

\begin{remark}\label{3-rem}
From Lemma \ref{4-thm}, we can read that $F(\cdot)\in D(G)$, with $G(x_1,\cdots,x_d)=\Pi_{i=1}^dG_1(x_i)$, and $F$ is asymptotically independent.
\end{remark}

Before we can prove Lemma~\ref{4-thm},  we present a modified version of Proposition~5.27 in Rensick~\cite[pp.296]{R1987}, which plays a key role in the proof of Lemma \ref{4-thm}.
\begin{lem}\label{4-lem-2}
Suppose that $\big\{X_n=(X_1^{(n)}\cdots,X_d^{(n)})',\;n\in\mathbb{N}\big\}$ are i.i.d.random vectors in $\R^d$ with the common joint continuous distribution $\tilde{F}(\cdot)$, and the marginal distributions $\tilde{F}_i(\cdot)$, $i=1,\cdots,d$. Moreover, we assume that $\tilde{F}_i(\cdot)$, $i=1,\cdots,d$ are both in the domain of attraction of some univariate extreme value distribution $\hat{G}_1(\cdot)$, i.e., there exist constants $a_i^{(n)}$ and $b_i^{(n)}$ such that
\beqnn
    \tilde{F}_i\Big(a_i^{(n)}x+b_i^{(n)}\Big)\to \hat{G}_1(x).
\eeqnn
Then, the following are equivalent:
\begin{itemize}
\item[(1)] $\tilde{F}$ is in the domain of attraction of a product measure, that is,
\beqnn
\tilde{F}^n\Big(a_i^{(n)}x_i+b_i^{(n)},i=1,\cdots,d\Big)\to \Pi_{i=1}^d \hat{G}_1\big(x_i\big);
\eeqnn
\item[(2)] For any $1\leq i<j\leq d$,  with $\lim_{x\to\infty}\tilde{F}_i(x)=1$
\beqlb\label{4-10}
\lim_{t\to\infty}\P\Big(X_i>t, X_j>t\Big)/\big(1-\tilde{F}_i(t)\big)\to 0.
\eeqlb
\end{itemize}
\end{lem}

By a slight modification of the proof to Proposition~5.27 in Rensick~\cite[pp.296]{R1987}, we can prove the above lemma, details of which are omitted here.

\medskip
Now, we are ready to prove Lemma~\ref{4-thm}.
\medskip

\noindent \underline{\proof of Lemma \ref{4-thm}:}
For the readers to follow the proof  below easily, we first  recall some notations we introduced in the previous sections. Recall that the reflection matrix $R$, the regulator process $L$ and the $d$-dimensional Brownian motion $X=(X_1,\cdots,X_d)'$ are the components in the definition of the SRBM given in (\ref{def1}), the time-change process $T$ is defined through (\ref{2-a3}),  and the sticky Brownian motion $Z$ is defined in (\ref{1-4}). Without loss of generality, we assume that $Z(0)={\bf 0}$. We mainly use the lemma \ref{4-lem-2} to prove this lemma.
Let
\beqnn
\hat{L}(t)=-[R^{-1}X(t)\wedge R^{-1}\mu t].
\eeqnn
Then, it follows from Konstantopoulos, Last and Lin \cite[Proposition 1]{KLL2004}  that  for any ${\bf \tilde{z}}=(\tilde{z}_1,\cdots,\tilde{z}_d)'\in\R^d_+$,
\beqlb
\P\{Z(t)\geq {\bf \tilde{z}}\}\leq \P\{\hat{Z}(t)\geq {\bf \tilde{z}}\},
\eeqlb
where $\hat{Z}(t)=\bar{Z}(T(t))$ with
\beqnn
\bar{Z}(t)=X(t)+R\hat{L}(t).
\eeqnn
It follows from  Kobayashi \cite[Lemma 2.7]{K2011} that
\beqnn
0<T^*:=\sup_{\omega}\{T(1,\omega)\}\leq 1,\;\text{a.s.}
\eeqnn
By the first change of variable formula (see, for example, Jacob \cite[Proposition 10.21]{J1979}), and the fact that $\hat{Z}(t)\geq 0$ and $\bar{Z}(t)\geq 0$ for all $t\in\R_+$, we have
\beqnn
\hat{Z}(1)=\int_0^{1} d \hat{Z}(s)=\int_0^{T(1)}d\bar{Z}(s)\leq \int_{0}^{T^*}d\bar{Z}(s)=\bar{Z}(T^*)\;\textrm{a.s.,}
\eeqnn
where the operations are performed component-wise. Hence, for any ${\bf \tilde{z}}=(\tilde{z}_1,\cdots,\tilde{z}_d)'\in\R^2_+$,
\beqlb\label{5-47}
\P\{Z(1)\geq {\bf \tilde{z}}\}\leq \P\{\bar{Z}(T^*)\geq {\bf \tilde{z}}\}.
\eeqlb
For convenience, let
\beqnn
\bar{F}({\bf \tilde{z}})=\P\big\{Z_1\geq \tilde{z}_1,\cdots, Z_d\geq \tilde{z}_d\}.
\eeqnn
We also note that
\beqlb\label{5-32}
\bar{F}({\bf\tilde{z}})=\lim_{t\to\infty}\P\{Z(t)\geq {\bf \tilde{z}}\}=\inf_{t\to\infty}\P\{Z(t)\geq {\bf \tilde{z}}\}\leq\P\{Z(1)\geq {\bf \tilde{z}}\}.
\eeqlb
From \eqref{5-47} and \eqref{5-32}, we get
\beqlb\label{5-4}
\bar{F}\big({\bf \tilde{z}}\big)\leq  \P\{\bar{Z}(T^*)\geq {\bf \tilde{z}}\}
\eeqlb
Below, we apply Lemma \ref{4-lem-2} to prove our result. Here, for convenience,  we assume that $i=1$ and $j=2$ in  Lemma \ref{4-lem-2}. Other cases can be discussed in the same fashion.  Furthermore,  for any ${\bf z}=(z_1,z_2)'\in\R_+^2$, let
\beqnn
\bar{F}_{12}({\bf z})=\P\{Z_1\geq z_1,\;Z_2\geq z_2\},
\eeqnn
and for a $d$-dimensional vector ${\bf {Y}}=(Y_1,\cdots,Y_d)'$,
\beqnn
{\bf Y}_{12}=(Y_1,Y_2)'.
\eeqnn
Therefore, from \eqref{5-4}, we get
\beqlb\label{4-11}
\bar{F}_{12}({\bf z})\leq \P\{\bar{Z}_1(T^*)\geq z_1, \bar{Z}_2(T^*)\geq z_2\}.
\eeqlb
On the other hand, for any ${\bf z}=(z_1, z_2)'\in\R^2_+$,
\beqlb\label{4-38}
\P\{Z_{12}(T^*)\geq {\bf z}\}&&\leq \P\{X_{12}(T^*)-\mu_{12} T^* \geq {\bf z}\}.
\eeqlb
It is obvious  that $X_{12}(T^*)-\mu_{12} T^*$ is a Gaussian vector with the correlation coefficient being less than 1.
\newline

From \eqref{4-38}, we have, for large enough $z\in\R_+$,
\beqlb\label{4-49}
\limsup_{z\to\infty}\frac{\bar{F}_{12}(z,z)}{\bar{F}_1(z)}\leq \limsup_{z\to\infty}\frac{ \P\{X_{12}(T^*)-\mu_{12} T^*\geq (z,z)'\}}{\bar{F}_1(z)}.
\eeqlb
At the same time, we know that if a bivariate Gaussian vector has a correlation coefficient strictly less than $1$, it is  asymptotically independent.  Hence,
\beqlb\label{4-38-a}
&&\limsup_{z\to\infty}\frac{\P\{X_1(T^*)-\mu_1 T^*\geq z, X_2(T^*)-\mu_2 T^*\geq z\}}{\P\{Z_1\geq z\}}\nonumber
\\
&&\hspace{2cm}=\limsup_{z\to\infty}\frac{\P\{X_1(T^*)-\mu_1 T^*\geq z, X_2(T^*)-\mu_2 T^*\geq z\}}{\P\{X_1(T^*)-\mu_1 T^*\geq z\}}\frac{\P\{X_1(T^*)-\mu_1 T^*\geq z\}}{\P\{Z_1\geq z\}}\nonumber
\\&&\hspace{2cm}\leq \limsup_{z\to\infty} \frac{\P\{X_1(T^*)-\mu_1 T^*\geq z, X_2(T^*)-\mu_2 T^*\geq z\}}{\P\{X_1(T^*)-\mu_1 T^*\geq z\}}\nonumber
\\&&\hspace{2cm}\leq \limsup_{z\to\infty} \frac{\P\{X_1(T^*)-\mu_1 T^*\geq z, X_2(T^*)-\mu_2 T^*\geq z\}}{\P\{X_1(T^*)-\mu_1 T^*\geq z\}}=0,
\eeqlb
where the first inequality is obtained by using the fact that
\beqnn
\P\{X_1(T^*)-\mu_1 T^*\geq z\}/\P\{Z_1\geq z\}\to 0,\;\text{as}\; z\to\infty.
\eeqnn
From above arguments, we obtain
\beqlb\label{4-39}
\lim_{z\to\infty}\frac{\bar{F}_{12}(z,z)}{\bar{F}_1(z)}=0.
\eeqlb
 From  \eqref{4-39} and Lemma \ref{4-lem-2}, the proof to the lemma follows.
\qed

From Lemma \ref{4-thm}, we can get the following result.

\begin{thm}\label{thm2}
For the multidimensional sticky Brownian motion $Z=(Z_1,\cdots,Z_d)'$,
\beqlb\label{5-36}
\P\big\{Z_1\geq z_i,\cdots, Z_d\geq z_d\big\}/\Big(\Pi_{i=1}^dg_i(z_i)\exp\{-\alpha_iz_i\}\Big)\to 1,
\eeqlb
as $(z_1,\cdots,z_d)'\to(\infty,\cdots,\infty)'$, where $g_i(\cdot)$ is given by \eqref{4-8}.
\end{thm}

\medskip

\noindent \underline{\proof of Theorem \ref{thm2}:}  To prove this theorem, we first introduce a transformation.
For the multivariate extreme value distribution $G(\cdot)$ defined in Remark \ref{3-rem},
\beqlb\label{5-10}
G^*(x_1,\cdots,x_d)=G\bigg(\Big(\frac{-1}{\log \big(G_1\big)}\Big)^{-1}\big(x_1),\cdots,\Big(\frac{-1}{\log \big(G_1\big)}\Big)^{-1}\big(x_d)\bigg).
\eeqlb
Then $G^*(\cdot)$ is the joint distribution function  with the common marginal Fr\'echnet distribution $\Phi(x)=\exp\{-x^{-1}\}$.
Furthermore,   for the stationary random vector ${\bf Z}=(Z_1,\cdots,Z_d)'$, define
\beqlb\label{5-12}
Y_i=\frac{1}{1-F_i(Z_i)}.
\eeqlb
Let $F^*(y_1,\cdots,y_d)$ be the joint distribution function of ${\bf Y}=(Y_1,\cdots,Y_d)'$. Then, it follows from Proposition 5.10 in Resnick \cite{R1987} and  Lemma~\ref{4-thm} that
\beqlb\label{5-6}
F^*(y_1,\cdots,y_d)\in D\big(G^*(y_1,\cdots,y_d)\big).
\eeqlb
By \eqref{5-6}, we have that for any ${\bf Y}=(y_1,\cdots,y_d)'\in\R_+^2$, as $n\to\infty$,
\beqlb\label{5-7}
(F^*(n{\bf Y}))^n\to G^*({\bf Y}) .
\eeqlb
It follows from \eqref{5-7} that
\beqnn
F^*(n{\bf Y})\sim \big(G^*({\bf Y})\big)^{\frac{1}{n}}.
\eeqnn
By a simple monotonicity argument, we can replace $n$ in the above equation by $t$.
Then we have, as $t\to\infty$,
\beqlb\label{5-3}
F^*(t{\bf Y})\sim \big(G^*({\bf Y})\big)^{\frac{1}{t}}.
\eeqlb
At the same time, by Lemma \ref{4-thm}, for any $y\in\R_+$,
\beqlb\label{5-28}
F^*_i(ty)\sim \big(G^*_1(y)\big)^{\frac{1}{t}},\;\text{for any}\;i=1,\cdots,d.
\eeqlb
Combining \eqref{5-3} and \eqref{5-28}, we get, as $t\to\infty$,
\beqlb\label{5-33}
F^*(t{\bf Y})\sim \Pi_{i=1}^dF_i^*(ty_i).
\eeqlb
It is obvious that for any $x\in\R_+$,
\beqlb\label{5-35}
\bar{F}^*_i(tx):=1-F_i^*(tx)\to 0\;\text{as}\; t\to \infty.
\eeqlb
Let $C(\bar{u}_1,\cdots,\bar{u}_d)$ be the copula of the random vector $(Y_1,\cdots,Y_d)'$, i.e.,
\beqlb\label{7-a27}
C\Big(F_1^*(z_1),\cdots,F_d^*(z_d)\Big)=F^*(z_1,\cdots,z_d).
\eeqlb
Furthermore, let $\hat{C}(\bar{u}_1,\cdots,\bar{u}_d)$ be the corresponding  survival copula of $C$. Then we have (see, for example, Schmitz \cite[Equation (2.46)]{SV2003} ):
\beqlb\label{7-a41}
&&\hat{C}(\bar{u}_1,\cdots,\bar{u}_d)=\sum_{i=1}^d \bar{u}_i+\sum_{1\leq i<j\leq n}C_{i,j}(1-\bar{u}_i,1-\bar{u}_i)-(n-1)\nonumber
\\&&\hspace{1cm}-\sum_{1\leq i<j<k\leq n}C_{i,j,k}(1-\bar{u}_i,1-\bar{u}_j,1-\bar{u}_k)+\cdots+(-1)^n C_{1,\cdots,d}(1-\bar{u}_1,\cdots,1-\bar{u}_n).
\eeqlb
For convenience, for any $(x_1,\cdots,x_d)'\in\R_+^d$, let $\bar{u}_i(t)=\bar{F}^*_i(tx_i)$. Hence, for any $t\in\R_+$,
\beqlb\label{7-a40}
\hat{C}\big(\bar{u}_1(t),\cdots,\bar{u}_d(t)\big)=\bar{F}^*(tx_1,\cdots,tx_d),
\\ C\big(1-\bar{u}_1(t),\cdots,1-\bar{u}_d(t)\big)=F^*(tx_1,\cdots,tx_d).\nonumber
\eeqlb
 Moreover, from \eqref{5-33}, we get, as $t\to\infty$,
\beqlb\label{7-a42}
C\big(1-\bar{u}_1(t),\cdots,1-\bar{u}_3(t)\big)\sim \Pi_{i=1}^d\big(1-\bar{u}_i(t)\big),
\eeqlb
and, for any $1\leq i_1<\cdots<i_k\leq 3$ with $k=2,\cdots,d$,
\beqlb\label{7-a43}
C_{i_1,\cdots,i_d}\big(1-\bar{u}_{i_1}(t),\cdots,1-\bar{u}_{i_k}(t)\big)\sim\Pi_{q=1} ^k\big(1-\bar{u}_{i_q}(t)\big).
\eeqlb
From \eqref{7-a41}, \eqref{7-a42} and \eqref{7-a43}, we can obtain that, as $t\to\infty$,
\beqlb\label{7-a24}
\hat{C}\big(\bar{u}_1(t),\cdots,\bar{u}_d(t)\big)\sim  \Pi_{i=1}^d\bar{u}_i(t),
\eeqlb
which, for any $(z_1,\cdots,z_d)'\in\R^d_+$, is equivalent to
\beqlb\label{5-44}
\lim_{t\to \infty}\frac{\bar{F}^*(tz_1,\cdots,tz_d)}{\Pi_{i=1}^d\bar{F}^*_i(tz_i)}=1.
\eeqlb
To prove our theorem, it suffices to show that
\beqlb\label{5-41}
\lim_{(z_1,\cdots,z_d)'\to(\infty,\cdots,\infty)'}\frac{\bar{F}^*(z_1,\cdots,z_d)}{\Pi_{i=1}^d\bar{F}^*_i(z_i)}=1.
\eeqlb
Note that
\beqlb\label{5-9}
\bar{F}^*(z_1,\cdots,z_d)=\P\big\{\bar{F}_1^*(Y_1)\geq \bar{F}_1^*(z_1), \cdots,\bar{F}_d^*(Y_d)\geq \bar{F}_d^*(z_d) \big\}.
\eeqlb
From \eqref{7-a27}, to prove \eqref{5-41}, we only need to show that
\beqlb\label{7-a28}
\lim_{(\bar{u}_1,\cdots,\bar{u}_d)'\to(0,\cdots,0)'\;\textrm{and}\; (\bar{u}_1,\cdots,\bar{u}_d)'\in I^d}\frac{\hat{C}(\bar{u}_1,\cdots,\bar{u}_d)}{\Pi_{i=1}^d\bar{u}_i}=1,
\eeqlb
where $I=[0,\;1]$.
We also recall that
\beqlb\label{7-a29}
\lim_{x\to 0}\frac{1-\exp\{-x\}}{x}=1.
\eeqlb
Hence, from \eqref{5-44} and  \eqref{7-a29}, we get, for any $(\bar{u}_1,\cdots,\bar{u}_d)'\in I^d$, that
\beqlb\label{7-7}
\lim_{t\to 0+}\frac{\hat{C}(t\bar{u}_1,\cdots,t\bar{u}_d)}{t^d\bar{u}_1\cdots \bar{u}_d}=1.
\eeqlb
Conversely, we note that the limit \eqref{7-a28} has the indeterminate form $\frac{0}{0}$. Hence, we would like to apply the multivariate L'h$\hat{o}$pital's rule (see Theorem 2.1 in \cite{L2012}) to prove it.  Without much effort, we can construct a multivariate differential function $\tilde{C}(\bar{u}_1,\cdots,\bar{u}_d)$, such that
\beqnn
\hat{C}(\bar{u}_1,\cdots,\bar{u}_d)=\tilde{C}(\bar{u}_1,\cdots,\bar{u}_d)\;\textrm{for all}\; (\bar{u}_1,\cdots,\bar{u}_d)'\in I^d,
\eeqnn
and
\beqnn
\tilde{C}(t\bar{u}_1,\cdots,t\bar{u}_d)\sim t^d \bar{u}_1\cdots \bar{u}_d,\;\textrm{as}\;t\to 0.
\eeqnn
Hence, it suffices to show that
\beqlb\label{7-8}
\lim_{(\bar{u}_1,\cdots,\bar{u}_d)'\to(0,\cdots,0)'\;\textrm{and}\; (\bar{u}_1,\cdots,\bar{u}_d)'\in I^d}\frac{\hat{C}(\bar{u}_1,\cdots,\bar{u}_d)}{\Pi_{i=1}^d \bar{u}_i}=\lim_{(\bar{u}_1,\cdots,\bar{u}_d)'\to(0,\cdots,0)'}\frac{\tilde{C}(\bar{u}_1,\cdots,\bar{u}_d)}{\Pi_{i=1}^d\bar{u}_i}=1.
\eeqlb
Near the origin $(0,\cdots,0)'$, the zero sets  of both $\tilde{C}(\bar{u}_1,\cdots,\bar{u}_d)$ and $\bar{u}_1\cdots \bar{u}_d$  consist of the hypersurfaces $\bar{u}_i=0$, $i=1,\cdots,d$.  By the multivariate L'h$\hat{o}$pital's rule  (see Theorem 2.1 in   \cite{L2012}), to prove \eqref{7-8}, it is enough to show that for each component $E_i$ of $\R^d\setminus\mathcal{C}$, where $\mathcal{C}=\cup_{i=1}^d\{\bar{u}_i=0\}$, we can find a vector $\vec{z}$, not tangent to $(0,\cdots,0)'$, such that, $D_{\vec{z}}(\Pi_{i=1}^d\bar{u}_i)\neq 0$ on $E_i$ and
\beqnn
\lim_{(\bar{u}_1,\cdots,\bar{u}_d)'\to (0,\cdots,0)'\; \textrm{and}\; (\bar{u}_1,\cdots,\bar{u}_d)'\in E_i} \frac{D_{\vec{z}}\tilde{C}(\bar{u}_1,\cdots,\bar{u}_d)}{D_{\vec{z}}(\bar{u}_1\cdots \bar{u}_d)}=1.
\eeqnn
For the component $E_1$ bounded by the hypersurfaces of $\mathcal{H}_i=\{(\bar{u}_1,\cdots,\bar{u}_d)':(\bar{u}_1,\cdots,\bar{u}_d)'\in\R_+^d\;\textrm{and}\; \bar{u}_i=0\}$,  for all $i=1,\cdots,d$, choose, say $\vec{z}=(1,\cdots,1)'$. Then, $z$ is not tangent to any hypersurfaces $u_i=0$, $i=1,\cdots,d$  at the point $(0,\cdots,0)'$.   Next, we take the limit along the direction $\vec{z}=(1,\cdots,1)'$. It follows from \eqref{7-a29} and  \eqref{7-7} that
\beqlb\label{7-a30}
\lim_{(\bar{u}_1,\cdots,\bar{u}_d)'\to (0,\cdots,0)'\;\textrm{and}\;(\bar{u}_1,\cdots,\bar{u}_d)'\in E_1} \frac{D_{\vec{z}}\tilde{C}(\bar{u}_1,\cdots,\bar{u}_d)}{D_{\vec{z}}(\bar{u}_1\cdots \bar{u}_d)}=1.
\eeqlb
Similar to \eqref{7-a30}, for any other components $E_i$, $i=2,\cdots,2^d$, we can find a vector $\vec{z}$ such that $z$ is not tangent to any hypersurfaces $\bar{u}_i=0$, $i=1,\cdots,d$  at the point $(0,\cdots,0)'$. Moreover, we have
\beqlb\label{7-a31}
\lim_{(\bar{u}_1,\cdots,\bar{u}_d)'\to (0,\cdots,0)'\;\textrm{and}\;(\bar{u}_1,\cdots,\bar{u}_d)'\in E_i} \frac{D_{\vec{z}}\tilde{C}(\bar{u}_1,\cdots,\bar{u}_d)}{D_{\vec{z}}(\bar{u}_1\cdots \bar{u}_d)}=1.
\eeqlb
From \eqref{7-8} to \eqref{7-a31} and   Theorem 2.1 in \cite{L2012},
\beqlb\label{7-a33}
\lim_{(\bar{u}_1,\cdots,\bar{u}_d)'\to (0,\cdots,0)'\;\textrm{and}\;(\bar{u}_1,\cdots,\bar{u}_d)'\in I^d} \frac{\hat{C}(\bar{u}_1,\cdots,\bar{u}_d)}{\bar{u}_1\cdots \bar{u}_d}=1.
\eeqlb

\noindent Finally, it follows from \eqref{5-12} that for any $(z_1,\cdots, z_d)'\in\R_+^d$,
\beqlb\label{5-13}
\P\{Z_1\geq z_1,\cdots,Z_d\geq z_d\}&&=\P\big\{F_1(Z_1)\geq F_1(z_1),\cdots,F_d(Z_d)\geq F_d(z_d)\big\}\nonumber
\\&&=\P\big\{Y_1\geq \frac{1}{1-F_1(z_1)},\cdots,Y_d\geq \frac{1}{1-F_d(z_d)}\big\}\nonumber
\\&&=F^*\Big(\frac{1}{\bar{F}_1(z_1)},\cdots,\frac{1}{\bar{F}_d(z_d)}\Big).
\eeqlb
Combining \eqref{5-41} and \eqref{5-13}, we get
\beqlb\label{5-15}
\P\{Z_1\geq z_1,\cdots,Z_d\geq z_d\}/\Bigg( \Pi_{i=1}^d\bar{F}^*_1\bigg(\frac{1}{\bar{F}_1(z_1)}\bigg)\Bigg)\to 1,\;\text{as}\; (z_1,\cdots,z_d)'\to(\infty,\cdots,\infty)'.
\eeqlb
 By \eqref{7-7} and \eqref{5-15}, we obtain
\beqlb\label{5-40}
\P\{Z_1\geq z_1,\cdots,Z_d\geq z_d\}/\Big( \Pi_{i=1}^d\bar{F}_i(z_i)\Big)\to 1,\;\text{as}\; (z_1,\cdots,z_d)'\to (\infty,\cdots,\infty)'.
\eeqlb
From the above arguments, the theorem is proved. \qed

\section{Exact Tail Asymptotic for Some Special Cases} \label{sec:4}

In the previous section, we studied the rough decay rate of multidimensional sticky Brownian motion. However, we are also interested in exact tail asymptotics for a multidimensional sticky Brownian motion.  In general, it is very difficult to get such results.  In this section, we study exact tail asymptotics for some special cases.

\noindent{\bf Example  4.1} ({\em Two-dimensional case}): We first consider the case for $d=2$. Dai and Zhao \cite{DZ2018} has proved that the marginal stationary distributions have the following exact tail asymptotics:
\beqlb\label{4-15}
\P\{Z_i\geq z_i\}\sim K_i z_{i}^{-\beta_i }\exp\{-\alpha_i z_i\}, \;i=1,2,
\eeqlb
where $\beta_i\in \{\frac{1}{2},\frac{3}{2},-1,0\}$ and $K_i$ is a non-zero constant.   By using the same method that was using in the proof to Theorem \ref{thm2} ,  we can show that, for a two-dimensional sticky Brownian motion, we have, as $(z_1,z_2)'\to(\infty,\infty)'$,
\beqnn
\P\big\{Z_1\geq z_1, Z_2\geq z_2\big\}/\Big(\Pi_{i=1}^2K_i z^{-\beta_i}_i\exp\{-\alpha_iz_i\}\Big)\to 1.
\eeqnn

\noindent{\bf Example 4.2} ({\em Skew symmetry case}):
We next consider the skew symmetric case.  In  R\'acz and Shkolnikov \cite[Theorem 5]{RS2015}, it was demonstrated that the stationary distribution for a multidimensional  sticky Brownian motion in a wedge  admits a separable form if the data satisfies  the conditions (14) and (15) in R\'acz and Shkolnikov \cite{RS2015}.  Following the proof to Theorem 5 in \cite{RS2015}, and noting the skew conditions for the SRBM on a  nonnegative orthant $\R_+^d$, (see, for example, \cite{ADH2001,CY2001,HW1987,HW1987b}), we can see that the sticky Brownian motion $Z$  has a separable  form if it satisfies the following skew symmetry condition:
\beqlb\label{4-16}
2\Sigma= R \Delta_{R}^{-1}\Delta_{\Sigma}+\Delta_{\Sigma}\Delta_{R^{-1}}R^{T},
\eeqlb
where $\Delta_A$ is the diagonal matrix with diagonal entries of a square matrix $A$.   In this case,  we can easily get that, as $(z_1,\cdots, z_d)'\to(\infty,\cdots,\infty)'$,

\beqnn
\P\big\{Z_1\geq z_1, \cdots, Z_d\geq z_d\big\}/\Big(\Pi_{i=1}^dK_i\exp\{-\alpha_iz_i\}\Big)\to 1.
\eeqnn

\noindent{\bf Example 4.3} ({\em Decomposability}): We first note that under the skew symmetry condition \eqref{4-16}, the stationary distribution of $Z$ could be obtained explicitly via the method used in R\'acz and Shkolnikov \cite[pp.1169-1170]{RS2015}. However, the condition \eqref{4-16} may be too strong, and is not satisfied in most cases.   At the same time, from the equation \eqref{2-a3},  we easily  see that when $u_i=0$, for all $i=1,\cdots,d$, the sticky Brownian motion $Z$ becomes an SRBM on a nonnegative orthant $\R_+^d$.   For the SRBM,  in applications, even if the condition \eqref{4-16} is not satisfied,  we can apply the product form based approximation to study the stationary distribution, see for example,\cite{KFLW2009}.  This product form based approximation may be improved by the decomposability in Dai, Miyazawa and Wu \cite{DMW2015}.   At the end of this work, we consider this special case, that is, $Z$ is an SRBM with the data $(\Sigma,\mu,R)$.  Let $\mathbb{J}:=\{1,\cdots,d\}$. We say that a pair $(\mathbb{K},\mathbb{L})$ is a partition of $\mathbb{J}$ if it satisfies $\mathbb{K}\cup\mathbb{L}=\mathbb{J}$  and  $\mathbb{K}\cup\mathbb{L}=\emptyset$. Recall that if the stationary distribution is the product of two marginal distributions associated with a partition $(\mathbb{K},\mathbb{L})$ of the set $\mathbb{J}$, then the stationary distribution is said to be decomposable with respect to $\mathbb{K}$ and $\mathbb{L}$.  Let $A^{(\mathbb{K},\mathbb{L})}$ be the $|\mathbb{K}|\times|\mathbb{L}|$ submatrix of a $d$-dimensional square matrix $A$ whose row and column indices are taken from $\mathbb{K}$ and $\mathbb{L}$, respectively, and $x^{\mathbb{K}}$ be the $\mathbb{K}$-dimensional vector with $x_i^{\mathbb{K}}=x_i$  for $i\in\mathbb{K}$, where  $x_i^{\mathbb{K}}$  is the $i$-th entry of   $x^{\mathbb{K}}$ . From Theorem 2 in Dai, Miyazawa and Wu \cite{DMW2015}, we get that if the covariance matrix $\Sigma$ and $R$ satisfy
\beqlb\label{4-17}
2\Sigma^{(\mathbb{K},\mathbb{K})}=R^{(\mathbb{K},\mathbb{K})} \Delta_{(R^{(\mathbb{K},\mathbb{K})})^{-1}} \Delta_{\Sigma^{(\mathbb{K},\mathbb{K})}}+\Delta_{\Sigma^{(\mathbb{K},\mathbb{K})}}\Delta_{(R^{(\mathbb{K},\mathbb{K})})^{-1}}(R^{(\mathbb{K},\mathbb{K})})^{T},
\eeqlb
\beqlb\label{4-18}
2\Sigma^{(\mathbb{L},\mathbb{K})}=R^{(\mathbb{L},\mathbb{K})} \Delta_{(R^{(\mathbb{K},\mathbb{K})})} \Delta_{(R^{(\mathbb{K},\mathbb{K})})^{-1}},
\eeqlb
and if the $|\mathbb{L}|$-dimensional $\Big(\Sigma^{(\mathbb{L},\mathbb{K})},\tilde{\mu}(\mathbb{L}), R^{(\mathbb{L},\mathbb{K})}\Big)$-SRBM has a stationary distribution, where $\tilde{\mu}(\mathbb{L})=\big(Q^{(\mathbb{L},\mathbb{L})})=\big(Q^{(\mathbb{L},\mathbb{L})}\big)^{-1}(Q\mu)^{\mathbb{L}}$ with $Q=R^{-1}$, then $Z^{\mathbb{K}}(0)$ and $Z^{\mathbb{L}}(0)$ are independent and $Z^{\mathbb{K}}(0)$ is of product form under $\pi$.  Hence, let $\mathbb{L}\subset \mathbb{J}$ with $|\mathbb{L}|=2$.  If the above conditions are satisfied, then we can get,  from Examples 4.1 and 4.2, that as
 $(z_1,\cdots, z_d)'\to(\infty,\cdots,\infty)'$,

\beqnn
\P\big\{Z_1\geq z_1, \cdots, Z_d\geq z_d\big\}/\Big(\Pi_{i=1}^dK_i z^{-\beta_i}_i\exp\{-\alpha_iz_i\}\Big)\to 1.
\eeqnn

\noindent\\[0.1mm]

\noindent{\bf Acknowledgments:} This research work was supported by the Fostering Project of Dominant Discipline and Talent Team of Shandong Province Higher Education Institutions, and the Natural Sciences and Engineering Research Council (NSERC) of Canada.

\end{document}